\documentclass[journal,twoside,web]{ieeecolor}
\usepackage{generic}
\usepackage{cite}
\usepackage{amsmath,amssymb,amsfonts}
\usepackage{algorithmic}
\usepackage{graphicx}
\usepackage{textcomp}
\usepackage{mathrsfs}
\usepackage{amstext}
\def\BibTeX{{\rm B\kern-.05em{\sc i\kern-.025em b}\kern-.08em
    T\kern-.1667em\lower.7ex\hbox{E}\kern-.125emX}}
\begin{document}
\title{Stochastic Adaptive Linear Quadratic Differential Games}
\author{Nian Liu and Lei Guo, \IEEEmembership{Fellow, IEEE}
\thanks{This paper was supported by the National Natural Science Foundational of China under Grant No.11688101.}
\thanks{The authors are with Institute of Systems Science, AMSS, Chinese
	Academy of Sciences, Beijing, 100190, China. (e-mails: liunian@amss.ac.cn,
	Lguo@amss.ac.cn).}}

\maketitle

\begin{abstract}
	 Game theory is playing more and more important roles in understanding complex systems and in investigating intelligent machines with various uncertainties. As a starting point, we consider the classical two-player zero-sum linear-quadratic stochastic differential games, but in contrast to most of the existing studies, the coefficient matrices of the systems are assumed to be unknown to both players, and consequently it is necessary to study adaptive strategies of the players, which may be termed as adaptive games and which has rarely been explored in the literature. In this paper, it will be shown that the adaptive strategies of both players can be constructed by the combined use of a weighted least squares (WLS) estimation algorithm, a random regularization method and a diminishing excitation method. Under almost the same structural conditions as those in the traditional known parameters case, we will show that the closed-loop adaptive game systems will be globally stable and asymptotically reaches the Nash equilibrium as the time tends to infinity. 

\end{abstract}


\section{Introduction}
\label{sec:introduction}
\IEEEPARstart{C}{omplex} systems are currently at the research frontiers of many fields in scientific technology, such as economic and social systems, biology and environmental systems, physical and engineering systems, and artificial intelligent systems. It is quite common that the components or subsystems of complex systems have game-like relationships, and the theory of differential games appears to be a useful tool in modeling and analyzing conflicts in the context of dynamical systems. 

 The differential game theory was firstly introduced in combat problems \cite{b8}, and has been applied in many fields (see, e.g.,\cite{b9,t1,t2,t3}). A great deal of research effort has been devoted to the area in the past half a century and much progress has been made (see, e.g., \cite{j2,j4,b1}). In particular, the linear-quadratic differential games, which are described by linear systems and quadratic payoff functions, have attracted a lot of attention. Bernhard \cite{p9} gave necessary and sufficient conditions for the existence of a saddle point for deterministic two-player zero-sum differential games on a finite time interval. Starr and Ho \cite{p10} extended the zero-sum differential games to the general case, i.e.,  the players wish to minimize different performance criteria and they discussed three types of solutions. Different state information patterns may give rise to different types of equilibrium, which have been discussed extensively in \cite{b1}. For examples, one of the classical result is feedback Nash equilibrium in linear quadratic differential games over an infinite time horizon, where the existence of Nash equilibrium is equivalent to the existence of solutions to a set of algebraic Riccati equations (see \cite{p11} and \cite{b1}). Another classical result is Nash equilibrium in open-loop linear quadratic Stackelberg differential games with sufficient condition for the existence of Nash equilibrium in terms of the existence of stabilizing solutions to three coupled algebraic Riccati equations (see \cite{b5}). Differential game theory may also be used to study $H^{\infty}$-optimal control problems \cite{p12}, since this problem is actually a minimax optimization problem. However, in all the above-mentioned works, the parameters in the related mathematical models are assumed to be known to the players.

When the parameters in stochastic dynamical systems are unknown, there are a great deal of researches in the area of adaptive control and much progress has been made over the past half a century. A basic method in the design of adaptive control is called the certainty equivalent principle, which consists of two steps: firstly to use the observed information to get an estimate of the unknown parameters at each time instant, then to construct or update the controller by taking the estimate as ``true" parameters at the same time. This design method is well-known to be quite powerful in dealing with dynamical systems with  possible large uncertainties (see \cite{j9,b6,b12}). However, since the closed-loop systems of adaptive control are usually described by a set of very complicated nonlinear stochastic dynamical equations, a rigorous theoretical investigation is well-known to be quite hard, even for linear uncertain stochastic systems. For example, one of the most well-known open problem in adaptive control has been to establish a rigorous convergence theory for the classical self-tuning regulators which are designed based on the least-squares estimates for linear stochastic systems with unknown parameters (see \cite{j10,l5,p1}). Another example is the adaptive Linear-Quadratic-Gaussian (LQG) control problem, where a key theoretical difficulty was how to guarantee the controllability of the online estimated model. This longstanding problem was reasonably resolved in the work \cite{p2} based on the self-convergence property of a class of weighted least squares and on a random regulation theory established in \cite{p3}, which turn out to be the fundamental bases for solving the adaptive game problems in the current paper. 

It goes without saying that uncertainties in the system structure, information and environment widely exist in dynamical games, and it is thus natural to consider adaptive game theory. To the best of our knowledge, only a little literature has been devoted to adaptive game theory, due to the complexity of the related theoretical investigation. For examples, Li and Guo\cite{p17} had considered a two-player zero-sum stochastic adaptive differential linear-quadratic game with state matrix to be known and stable. Yuan and Guo\cite {j7} investigated adaptive strategies for a zero-sum game described by an input-output stochastic model with known high gain parameters for both players. In a related but somewhat different context, Li, Liu and Wang\cite{p16} studied the linear-quadratic $H^{\infty}$-optimal control problem with unknown system parameters in a framework of zero-sum differential game, by using an integral reinforcement learning-based adaptive dynamic programming procedure.

In this paper, we consider the problem of stochastic adaptive linear-quadratic zero-sum differential games, using only the natural assumption of controllability and antianalytic property of the underlying uncertain stochastic systems. Firstly, the WLS scheme is used to obtain an online convergent estimate with no assumptions on the system signals (the self-convergence property). Secondly, the WLS estimate is modified by a random regularization procedure that was introduced in \cite{p3} to obtain a uniformly controllable family of estimated models. Thirdly, adaptive strategies that are globally stabilizing and convergent are designed for both players by using a diminishing excitation technique. It is showed that under the above adaptive strategies, the closed-loop game systems will reach the Nash equilibrium asymptotically, and in the meantime the estimates converge to the true parameters globally and almost surely.

The remainder of the paper is organized as follows: In Section II,  the design procedure of the adaptive strategies is provided and the main results on global stability, convergence of the estimate and Nash equilibrium of the closed-loop game systems are presented. Section III gives the proofs of above theorems and Section IV concludes this paper. Some proofs of lemmas are provided in the Appendix.

\section{PROBLEM FORMULATION AND MAIN RESULTS}
\subsection{Problem Formulation}
Consider the following basic stochastic linear quadratic zero-sum differential game: 
\begin{equation}
	dx(t)=(Ax(t)+B_1u_1(t)+B_2u_2(t))dt+Ddw(t),
\end{equation}
where $x(t)\in \mathbb{R}^n$ is the state with the initial state $x(0)=x_0$, $u_i\in\mathbb{R}^{m_i}$ is the strategy of Player $i$ $(i=1,2)$, $(w(t),\mathscr{F}_t;t\ge0)$ is a $\mathbb{R}^p$-valued standard Wiener process. $A\in \mathbb{R}^{n\times n},B_i\in\mathbb{R}^{n\times m_i},D\in \mathbb{R}^{n\times p}$ are system matrices.

The payoff function is :
\begin{equation}
	J(u_1,u_2)=\limsup\limits_{T\to \infty}\frac{1}{T}\int_{0}^{T}(x^{\mathsf{T}}Qx+u_1^{\mathsf{T}}R_1u_1-u_2^{\mathsf{T}}R_2u_2)dt,
\end{equation}
where $Q=Q^{\mathsf{T}},R_1=R_1^{\mathsf{T}}>0,R_2=R_2^{\mathsf{T}}>0$ are given weighting matrices. The objective of Player 1 is to minimize the payoff function, while Play 2 wants to maximize it. 

It is well-known that the information structure and the timing of actions play a crucial role in games (see \cite{b1}). The basic information structures of the above zero-sum game is as follows:

1) The system matrices $A,B_1,B_2,D$ are unknown to both players. But, the weighting matrices $\{Q,R_1,R_2\}$ and the state $x(t)$ with the initial state $x(0)$ are ``common" knowledge (see \cite{b2}).

2) The strategies $u_1(t)$ and $u_2(t)$ are of feedback patterns, i.e., $u_1$ and $u_2$ are adapted to $\{\mathscr{F}^x_t;t\ge0\}$ where $\mathscr{F}^x_t$ is a filtration generated by the  state process $x(t)$, i.e., $\mathscr{F}^x_t\triangleq \sigma(x(s);0\le s\le t)$.

\emph{Definition 1:} The strategy pair $(u_1(t),u_2(t))$ is said to be admissible if it is adapted to $\{\mathscr{F}^x_t;t\ge0\}$ and under which the system (1) is globally stable in the sense that for any initial state $x(0)$,
\begin{equation}\notag
	\limsup\limits_{T\to \infty}\frac{1}{T}\int_{0}^{T}(|x(s)|^2+|u_1(s)|^2+|u_2(s)|^2)ds<\infty\quad a.s.
\end{equation}

\emph{Remark 1:} The stability property is a natural requirement for the game, which is imposed to ensure the finiteness of the infinite-horizon cost integral that will be considered in the paper. Obviously, the stability depends on both strategy spaces of the players, which can be justified from the supposition that both players have a first priority in stabilizing the system as in the non-adaptive deterministic case \cite{b5}.

\emph{Definition 2:}  For the zero-sum linear-quadratic differential game with both players in the feedback pattern, a pair of admissible strategy $(u_1^*,u_2^*)$ is called a feedback Nash equilibrium if it satisfies
$$J(u_1^*,u_2)\le J(u_1^*,u_2^*)\le J(u_1,u_2^*),$$
for any admissible pairs $(u_1^*,u_2)$ and $(u_1,u_2^*)$.

It is well-known that if the algebraic Riccati equation  (ARE)
\begin{equation}\label{ra}
	A^{\mathsf{T}}P+PA+Q-PB_1R_1^{-1}B_1^{\mathsf{T}}P+PB_2R_2^{-1}B_2^{\mathsf{T}}P=0,
\end{equation}
admits a symmetric solution $P$ such that $A-(B_1R_1^{-1}B_1^{\mathsf{T}}-B_2R_2^{-1}B_2^{\mathsf{T}})P$ is stable, then the following pair of strategy is a feedback  Nash equilibrium:
\begin{equation}\label{h1}
	u_1(t)=-R_1^{-1}B_1^{\mathsf{T}}Px(t),
\end{equation}
\begin{equation}\label{h2}
u_2(t)=R_2^{-1}B_2^{\mathsf{T}}Px(t).
\end{equation}
It is worth mentioning that such solution $P$ is called the stabilizing solution to ARE (\ref{ra}) and it is unique (see \cite{b5}).

\emph{Definition 3:} A matrix function $G(s)$ is called antianalytic facotorizable if there exist matrix functions $\Xi(s)$ and $\Omega(s)$ such that
$$G(j\omega)=\Xi(-j\omega)\Omega(j\omega),$$
where $\Xi,\Xi^{-1},\Omega,\Omega^{-1}$ are all stable proper rational matrix functions. More details  concerning antianalytic  factorizations can be found in \cite{p13} and \cite{j8}.

Now, we introduce a class of matrices defined by
$$\mathscr{F}\triangleq\bigg\{L\triangleq\left(\begin{array}{c}
	L_1\\
	L_2\\
\end{array}\right)|A+B_1L_1+B_2L_2\text{ is stable}\bigg\},$$
and we need the following notations:
\begin{equation}
	B=(B_1,B_2),\quad R=diag(R_1,-R_2),
\end{equation}
\begin{equation}\label{r3}
	G(s)=R+B^{\mathsf{T}}(-sI-A^{\mathsf{T}})^{-1}Q(sI-A)^{-1}B,
\end{equation} 
\begin{equation}
	N_L(s)=I+L(sI-A-BL)^{-1}B,
\end{equation}
\begin{equation}
	\tilde{G}_L(s)=N_L^*(s)G(s)N_L(s),
\end{equation}
where $N^*$ represent the conjugate transpose of $N$.

\emph{Definition 4:} Assume $(A, B)$ is stabilizable. We will say that $G(s)$ is antianalytic prefactorizable (see \cite{p13})  if there exists $L\in\mathscr{F}$ such that $\tilde{G}_L(s)$ is antianalytic  factorizable.

\emph{Proposition 1\cite{p13}:} The algebraic Riccati equation (\ref{ra}) has a stabilizing solution and $R$ is nonsingular if and only if the pair $(A,B)$ is stabilizable and the matrix function $G(s)$ is antianalytic perfactorizable.

Then, the following basic assumptions are made.

A1) the pair $(A,B)$ is controllable.

A2) the matrix function $G(s)$ defined by (\ref{r3}) is antianalytic perfactorizable.

\emph{Remark 2:} Assumption A1) and A2) ensure that the ARE (\ref{ra}) has a stabilizing solution. Unfortunately, the corresponding strategies (\ref{h1}) and (\ref{h2}) are not implementable, because the system matrices $(A,B)$ are unknown to the players. To solve this problem, it is natural to use adaptive methods based on online estimation of the system matrices. The controllability assumption A1) also makes it possible to transform excitation properties from the input to the state signals, necessary for the convergence of the parameter estimation. These are the contents of the next subsections.

\subsection{The WLS Estimation}
To describe the estimation problem in the standard form, we introduce the following notations:
\begin{equation}\label{k1}
	\theta^{\mathsf{T}}=[A,B_1,B_2],
\end{equation}
and
$$\varphi(t)=[x^{\mathsf{T}}(t),u_1^{\mathsf{T}}(t),u_2^{\mathsf{T}}(t)]^{\mathsf{T}},$$
and rewrite the system (1) into the following compact form:
\begin{equation}
	dx(t)=\theta^{\mathsf{T}}\varphi(t)dt+Ddw(t).
\end{equation}

Now the continuous-time weighted least-square (WLS) estimates, $(\theta(t),t\ge0)$, are given by \cite{p2}
\begin{equation}\label{WLS1}
	d\theta(t)=a(t)Q(t)\varphi(t)[dx^{\mathsf{T}}(t)-\theta^{\mathsf{T}}\varphi(t)dt],
\end{equation}
\begin{equation}\label{WLS2}
	dQ(t)=-a(t)Q(t)\varphi(t)\varphi^{\mathsf{T}}(t)Q(t)dt,
\end{equation}
where $Q(0)>0$ and $\theta^{\mathsf{T}}(0)=[A(0),B_1(0),B_2(0)]$ are arbitrary deterministic values such that $(A(0),B(0))$ is controllable with $B(0)=[B_1(0),B_2(0)],$

\begin{equation}
	a(t)=\frac{1}{f(r(t))},
\end{equation}
\begin{equation}\label{r}
	r(t)=\|Q^{-1}(0)\|+\int_{0}^{t}|\varphi(s)|^2ds,
\end{equation}
and $f\in\mathbb{F}$ with 
	$$\mathbb{F}=\{f|f:\mathbb{R}_+\to\mathbb{R}_+,f \text{ is slowly increasing and }$$
\begin{equation}
	\int_{c}^{\infty}\frac{dx}{xf(x)}<\infty \text{ for some } c\ge0\}
\end{equation}
where a function is called slowly increasing if it is increasing
and satisfies $f\ge1$ and $f(x^2)=O(f(x))$ as $x\to \infty$ (see \cite{p3}).

\emph{Remark 3: }Some typical functions in $\mathbb{F}$ are, for example, 
$$f(x)=log^{(1+\delta)}x, (logx)(loglogx)^{(1+\delta)},\cdots \quad (\delta>0).$$

 \subsection{Regularization}
 To construct the adaptive version of the ARE (\ref{ra}) by using the estimates $(A(t),B(t))$ with guaranteed solvability, one needs at least the stabilizability of the estimates $(A(t),B(t))$ which may not be provided by the above WLS algorithm. To solve this problem, We resort to the regularization method introduced in \cite{p3} to modify the WLS estimates to ensure their uniform controllablility. We first introduce the following definition \cite{p2}:
 
\emph{Definition 4:} A family of system matrices $(A(t),B(t); t \ge 0)$ is said to be uniformly controllable if there is a constant $c > 0$ such that
$$\sum_{i=0}^{n-1}A^i(t)B(t)B^{\mathsf{T}}(t)A^{i{\mathsf{T}}}(t)\ge cI,$$
for all $t \in [0,\infty)$, where $A(t)\in R^{n\times n},B(t) \in R^{n\times m}$.

\emph{Lemma 1: }Let $(\theta(t),t \ge 0)$ be defined by (\ref{WLS1}) and (\ref{WLS2}). Then the
following properties are satisfied:
\begin{align}
	&(1) \sup\limits_{t\ge0}\|Q^{-1}(t)\tilde{\theta}(t)\|^2<\infty\quad a.s.\notag\\
	&(2)\int_{0}^{\infty}a(t)\|\tilde{\theta}^{\mathsf{T}}(t)\varphi(t)\|^2dt<\infty\quad a.s.\notag\\
	&(3)\lim\limits_{t\to\infty}\theta(t)=\bar{\theta}\quad a.s.\notag
\end{align}
where $\tilde{\theta}=\theta(t)-\theta$ and $\bar{\theta}$ is a matrix-valued random variable.
The proof of the above lemma can be found in \cite{p2}.

By Lemma 1(3), it is known that $\theta(t)$ converges to a certain random matrix $\bar{\theta}$ which may not be the true parameter matrix $\theta$ and naturally, the controllability of the estimate models may not be guaranteed. To solve this problem, we observe that by Lemma 1(1), the matrix sequence $\{Q^{-\frac{1}{2}}(t)(\theta-\theta(t)),\,\,t\ge0\}$ is bounded. This inspired the following modification for the WLS estimate (see \cite{p2} and \cite{p3}):
$$\theta(t,\beta(t))=\theta(t)-Q^{\frac{1}{2}}(t)\beta(t),$$
where $\beta(t)\in\mathbb{R}^{(n+m_1+m_2)\times n} $ is a sequence of bounded matrices to be defined shortly. For simplicity, we denote
$$\theta^{\mathsf{T}}(t,\beta(t))=[\bar{A}(t),\bar{B}(t)], \quad \bar{B}(t)=[\bar{B}_1(t),\bar{B}_2(t)].$$

To guarantee the uniform controllability of $(\bar{A}(t),\bar{B}(t))$, we need only to select the sequence $\{\beta(t),t\ge0\}$ to guarantee the uniform positivity of $Y(t)$ defined by
$$Y(t,\beta(t))=\det\big(\sum_{i=0}^{n-1}\bar{A}^i(t)\bar{B}(t)\bar{B}^{\mathsf{T}}(t)\bar{A}^{i{\mathsf{T}}}(t)\big).$$

For this purpose, let $\{\eta_k\in\mathbb{R}^{(n+m_1+m_2)\times n},k\in \mathbb{N}\}$ be a sequence of independent random variables which are uniformly distributed in the unit ball for a norm of the matrices and are also independent of $(w(t),t\ge0)$. The procedure of choosing $\beta$ is recursively given by the following:

\begin{align}
	&\beta_0=0,\notag \\
	&\beta_k=\begin{cases}
		\eta_k,&\text{if } Y(k,\eta_k)\ge(1+\gamma)Y(k,\beta_{k-1})\\
		\beta_{k-1},&\text{otherwise}
	\end{cases}
\end{align}
where $\gamma\in(0,\sqrt{2}-1)$ is a fixed constant. Thus, a sequence of
regularized estimates $(\bar{\theta}_k, k\in\mathbb{N})$ can be defined by 
\begin{equation}\label{key1}
	\bar{\theta}_k=\theta(k)-Q^{-1/2}(k)\beta_k.
\end{equation}

Finally, the continuous-time estimates used for the design of adaptive strategies can be defined piecewise as follows:
\begin{equation}\label{key2}
	\hat{\theta}(t)=\bar{\theta}_k,
\end{equation}
for any $t \in (k, k +1]$ and for all $k\in\mathbb{N}$.

\subsection{The Main Result}
For simplicity, we rewrite the estimates given by (\ref{key2}) as
$$\hat{\theta}^{\mathsf{T}}(t)=[A(t),B_1(t),B_2(t)].$$

For any $k\in\mathbb{N}$, It is well-known that the following ARE:
$$A^{\mathsf{T}}(k)P(k)+P(k)A(k)+Q-P(k)B_1(k)R_1^{-1}B_1^{\mathsf{T}}(k)P(k)$$
\begin{equation}\label{r7}
	+P(k)B_2(k)R_2^{-1}B_2^{\mathsf{T}}(k)P(k)=0,
\end{equation}
has at most one Hermitian matrix solution $P(k)$ such that
$$A_{cl}\big(P(k)\big)\triangleq A(k)-B_1(k)R_1^{-1}B_1^{\mathsf{T}}(k)P(k)$$
\begin{equation}\label{a17}
	+B_2(k)R_2^{-1}B_2^{\mathsf{T}}(k)P(k)
\end{equation}
is stable (see \cite{b14}). We can rewrite such $P(k)$ as
\begin{equation}\label{a18}
	P(k)=P^1(k)+P^2(k)j
\end{equation}
where $P^1(k)$ is a real symmetric matrix and $P^2(k)$ is a real skew-symmetric matrix and $j^2=-1$.

Then, we construct the desired strategy pair by considering two cases separately. Case (i): If the ARE (\ref{r7}) has a Hermitian matrix solution $P(k)$ such that both $A_{cl}\big(P(k)\big)$ and $A_{cl}\big(P^1(k)\big)$ are stable, then Players 1 and 2 can use following strategies respectively:
\begin{equation}\label{key4}
	u_1(t)=-R_1^{-1}B_1^{\mathsf{T}}(k)P^1(k)x(t),\text{ for }t\in(k, k +1],
\end{equation}
\begin{equation}\label{key5}
	u_2(t)=R_2^{-1}B_2^{\mathsf{T}}(k)P^1(k)x(t),\text{ for }t\in(k, k +1].
\end{equation}
Case (ii): If the ARE (\ref{r7}) does not admit any Hermitian matrix solution $P(k)$ such that both $A_{cl}\big(P(k)\big)$ and $A_{cl}\big(P^1(k)\big)$ are stable, then we choose the following strategy pair for $t\in(k, k +1]$:
\begin{equation}\label{ra1}
	(u_1^{\mathsf{T}}(t),u_2^{\mathsf{T}}(t))^{\mathsf{T}}=-B(k)^{\mathsf{T}}W_k^{-1}(0,T_0)x(t),
\end{equation}
where  $W_k(0,T_0)=\int_{0}^{T_0}e^{-A(k)\tau}B(k)(e^{-A(k)\tau}B(k))^{\mathsf{T}}d\tau$, $T_0>0$ is an arbitrary constant and $B(k)=[B_1(k),B_2(k)]$.

\emph{Remark 4: }The strategy (\ref{ra1}) is chosen to stabilize the system (1) whenever the ARE (\ref{r7}) does not have a desired solution. Once the adaptive game system is globally stabilized, the diminishing excitation method \cite{b6} may be used to guarantee the strong consistency of the WLS, and consequently the strategy pair $(u_1,u_2)$ will have the form (\ref{key4}) and (\ref{key5}) when $t$ is large enough, which will be discussed in the following.

Following the ideas as used in \cite{p2}, \cite{p3} and \cite{b6}, we add some diminishing excitation to the strategies (\ref{key4}), (\ref{key5}) and (\ref{ra1}), i.e., for $t\in(k,k+1]$
\begin{equation}\label{a1}
	u^*_1(t)=L_1(k)x(t)+\gamma^{(1)}_k(v_1(t)-v_1(k)),
\end{equation}
\begin{equation}\label{a2}
	u^*_2(t)=L_2(k)x(t)+\gamma^{(2)}_k(v_2(t)-v_2(k)),
\end{equation}
where $[L_1(t),L_2(t)]=[-R_1^{-1}B_1^{\mathsf{T}}(k)P^1(k),R_2^{-1}B_2^{\mathsf{T}}(k)P^1(k)]$ in Case (i), and $[L^{\mathsf{T}}_1(t),L^{\mathsf{T}}_2(t)]^{\mathsf{T}}=-B(t)^{\mathsf{T}}W_k^{-1}(0,T_0)$ in Case (ii). The sequences $\{\gamma^{(i)}_k, k\in\mathbb{N}\}$, $i=1,2$ can be any
sequences satisfying the following:
$$\lim\limits_{k\to \infty}\gamma^{(1)}_k=0,\quad \lim\limits_{k\to \infty}\gamma^{(2)}_k=0,$$ 
$$\min\{\gamma^{(1)}_k,\gamma^{(2)}_k\}\ge\gamma_k=(\frac{\log k}{\sqrt{k}})^{\frac{1}{2}},$$
for simplicity, we choose $\gamma^{(1)}_k=\gamma^{(2)}_k=\gamma_k$ in the rest of the paper. The processes $(v_1(t),t \ge 0)$ and $(v_2(t),t \ge 0)$ are chosen as sequences of independent standard Wiener processes that are independent of $(w(t),t\ge0)$ and $(\eta_k,k\in\mathbb{N})$.

The proofs of the following theorems will be given in the next section.

\emph{Theorem 1: }Let A1) be satisfied. Then the solution $(x(t),t \ge 0)$ of the system (1) under the
adaptive strategies (\ref{a1}) and (\ref{a2}) of the players, is globally stable in the sense that for any initial state $x(0)$,
\begin{equation}
	\limsup\limits_{T\to \infty}\frac{1}{T}\int_{0}^{T}|x(s)|^2ds<\infty\quad a.s.
\end{equation}

\emph{Remark 5: } By Theorem 1, it is not difficult to see that the strategy pair $(u_1^*,u_2^*)$ is admissible. It is worth mentioning that the excitation does not influence the stability of the system, but will provide sufficient excitation needed for strong consistence of the WLS estimates.

\emph{Theorem 2: }Let $(\hat{\theta}(t),t\ge0)$ be the family of estimates
given by (\ref{key2}) and the players use the above strategies (\ref{a1}) and (\ref{a2}) in the system (1). If A1) is satisfied, then
\begin{equation}
	\lim\limits_{t\to\infty}\hat{\theta}(t)=\theta\quad a.s.
\end{equation}
where $\theta$ is the true system parameter define by (\ref{k1}).

\emph{Theorem 3: }Consider the stochastic game system (1) with the payoff function (2). If A1) and A2) are satisfied , then the above adaptive strategies (\ref{a1}) and (\ref{a2}) constitute a feedback Nash equilibrium. Moreover, we have
\begin{equation}
	J(u_1^*,u_2^*)=tr(D^{\mathsf{T}}PD)\quad a.s.
\end{equation}
where $P$ is the stabilizing solution of the ARE (\ref{ra}).

\section{Proof}

We first present a basic lemma on WLS, which is similar to Lemma 2 of \cite{p2}.

\emph{Lemma 2:} Let A1) be satisfied for the game system (1) with the payoff function (2). Then the family of regularized WLS estimates $(\hat{\theta}(t),t\ge0)$ defined by (\ref{WLS1})-(\ref{key2}) has the following properties.
\par 
(1) Self-convergence, that is, $\hat{\theta}(t)$ converges a.s. to a finite random matrix as $t\to\infty$.
\par 
(2) The family $(A(t),B(t);t\ge0)$ is uniformly controllable where $B(t)=[B_1(t),B_2(t)]$ and $[A(t),B(t)]=\hat{\theta}^{\mathsf{T}}(t)$.
\par 
(3) Semi-consistency, that is, as $t\to \infty$
$$\int_{0}^{t}|(\hat{\theta}(s)-\theta)^{\mathsf{T}}\varphi(s)|^2ds=o(r(t))+O(1)$$ 
where $r(t)$ is defined by (\ref{r}).
\vskip0.2cm

We also need the following lemmas whose proofs are given in the Appendix.
\vskip0.2cm
\emph{Lemma 3: }If the system $\dot{x}=Ax+Bu$ is controllable, then under the control $u=-B^{\mathsf{T}}W^{-1}(0,T_0)x(t)$, the closed-loop system is stable, where $W(0,T_0)=\int_{0}^{T_0}e^{-At}B(e^{-At}B)^{\mathsf{T}}dt$ and $T_0>0$ is any given constant.
\vskip0.2cm

\emph{Lemma 4: }For the processes $(v_1(t),t \ge 0)$ and $(v_2(t),t \ge 0)$ and the sequences $\{\gamma_k, k\in\mathbb{N}\}$ defined in (\ref{a1}) and (\ref{a2}), we have
	$$
	\limsup\limits_{N\to \infty}\frac{1}{N}\sum_{k=1}^{N}\int_{k}^{k+1}\gamma_k^2|v_1(t)-v_1(k)|^2dt=0\quad a.s.
	$$
	$$
	\limsup\limits_{N\to \infty}\frac{1}{N}\sum_{k=1}^{N}\int_{k}^{k+1}\gamma_k^2|v_2(t)-v_2(k)|^2dt=0\quad a.s.
	$$
\vskip0.2cm

\emph{Lemma 5: }Let $(x(t), t\ge0)$ be the solution of (1) with the adaptive strategies (\ref{a1}) and  (\ref{a2}). Then 
$$\frac{1}{N}\sum\limits_{k=0}^{N}|x(k)|^2=O(1)+o(\frac{r(N)}{N})\quad a.s.$$

The following lemma on continuity of the stabilizing solution of the ARE can be found in \cite{b13}.
\vskip0.2cm
\emph{Lemma 6: }With each triple $(A,S,Q)\in(\mathbb{C}^{n\times n})^3$ satisfying $S=S^{\mathsf{T}}$ and $Q=Q^{\mathsf{T}}$ to be piecewise continuous and locally bounded matrix functions
$$A,S,Q:\mathbb{D}\to\mathbb{C}^{n\times n}$$
on some interval $\mathbb{D}\in\mathbb{R}$, we associate the matrix function
$$E=\left(\begin{array}{cc}
	Q &  A^{\mathsf{T}}\\ 
	A & -S \\ 
\end{array}\right):\mathbb{D}\to\mathbb{C}^{2n\times 2n},$$
and the ARE
\begin{equation}\label{r8}
	A^{\mathsf{T}}P+PA+Q-PSP=0.
\end{equation}
Assume that for some $E=E_0$ there is a stabilizing solution $P_0$ for which (\ref{r8}) is fulfilled. Then there exists $r(E_0)>0$ such that for $E$ ranging $||E-E_0||<r(E_0)$, there is a unique analytic function $E\mapsto P(E)$ such
that $A-SP(E)$ is stable and $P(E)$ is a Hermitian matrix solution to the equation (\ref{r8}) satisfying $P(E_0)=P_0$.

\subsection{Proof of Theorem 1:}
By Lemma 2, there are random matrices $A(\infty),B_1(\infty)$ and $B_2(\infty)$ such that
$$\lim\limits_{t\to \infty}A(t)=A(\infty) \quad  a.s.$$
$$\lim\limits_{t\to \infty}B_1(t)=B_1(\infty) \quad a.s.$$
$$\lim\limits_{t\to \infty}B_2(t)=B_2(\infty) \quad a.s.$$
and that $(A(\infty),B(\infty))$ is controllable a.s. where $B(\infty)=[B_1(\infty),B_2(\infty) ]$.
\par 
Now, we verify the convergence of $L_1(t)$ and $L_2(t)$ by considering two cases separately. (i) We first consider the case where the ARE
$$A^{\mathsf{T}}(\infty)P(\infty)+P(\infty)A(\infty)-P(\infty)B_1(\infty)R_1^{-1}B_1^{\mathsf{T}}(\infty)P(\infty)$$
\begin{equation}\label{r2}
	+Q+P(\infty)B_2(\infty)R_2^{-1}B_2^{\mathsf{T}}(\infty)P(\infty)=0,
\end{equation}
 have a Hermitian matrix solution $P(\infty)$ such that both $A_{cl}\big(P(\infty)\big)$ and $A_{cl}\big(P^1(\infty)\big)$ defined by (\ref{a17})(\ref{a18}) are stable. By Lemma 6, there are solutions $P(k)$ for all large enough k such that 
$$\lim\limits_{t\to \infty}P(k)=P(\infty) \quad  a.s. $$
and both $A_{cl}\big(P(k)\big)$ and $A_{cl}\big(P^1(k)\big)$ are stable, which means that when $k$ is large enough, $L_1(k)$ and $L_2(k)$ will take the form 
$$L_1(k)=-R_1^{-1}B_1^{\mathsf{T}}(k)P^1(k),\quad L_2(k)=R_2^{-1}B_2^{\mathsf{T}}(k)P^1(k)$$ almost surely.
(ii) In the case where the algebraic Riccati equation (\ref{r2}) does not admit any Hermitian matrix solution $P(\infty)$ such that both $A_{cl}\big(P(\infty)\big)$ and $A_{cl}\big(P^1(\infty)\big)$ are stable, by the definition of $[L_1(t),L_2(t)]$ in (\ref{a1}) and (\ref{a2}), we will have for $k$ is large enough $$[L^{\mathsf{T}}_1(t),L^{\mathsf{T}}_2(t)]^{\mathsf{T}}=-B(k)^{\mathsf{T}}W_k^{-1}(0,T_0),\text{ for } t\in(k,k+1],$$  
which is also a continuous matrix function of the parameters $(A(t),B(t))$. 

To summarize, in both cases we have
 $$\lim\limits_{t\to \infty}L_1(t)=L_1(\infty) \quad  a.s.$$
 $$\lim\limits_{t\to \infty}L_2(t)=L_2(\infty) \quad  a.s.$$
for some random matrices $L_1(\infty)$ and $L_2(\infty)$.

For simplicity of the remaining descriptions, we denote
\begin{equation}\label{l2}
	\Phi(t)=A(t)+B_1(t)L_1(t)+B_2(t)L_2(t),
\end{equation}
then we have
\begin{equation}\label{l3}
	\lim\limits_{t\to \infty}\Phi(t)=\Phi(\infty) \quad  a.s.
\end{equation}
where
$$\Phi(\infty)=A(\infty)+B_1(\infty)L_1(\infty)+B_2(\infty)L_2(\infty).$$

By Lemmas 2, 3, 6 and the above analysis, we know that $\Phi(t)$ is stable for any $t\ge0$ and converges to a stable matrix a.s. Hence, there exist some uniformly bounded positive definite matrices $K(t)$ such that
\begin{equation}\label{key7}
	\Phi^{\mathsf{T}}(t)K(t)+K(t)\Phi(t)=-I.
\end{equation}

Substituting (\ref{a1}) and (\ref{a2}) into the system (1), we have
\begin{align}\label{s}
dx(t)&=(Ax(t)+B_1u^*_1(t)+B_2u^*_2(t))dt+Ddw(t)\notag\\
&=(\Phi(t)x(t)+\delta(t)+\gamma_k(v(t)-v(k)))dt+Ddw(t),\notag
\end{align}
where $\delta(t)=(\theta-\hat{\theta}(t))^{\mathsf{T}}\varphi(t)$, $v(t)=B_1v_1(t)+B_2v_2(t)$ and $t\in(k,k+1]$, $k\in\mathbb{N}$.
By Lemma 2(3), we have
\begin{equation}\label{z1}
	\int_{0}^{t}|\delta(s)|^2ds=o(r(t))+O(1).
\end{equation}
Then, applying the Ito's formula to $\langle K(t)x(t),x(t)\rangle$ where $\langle,\rangle$ represents inner product, and noting that $K(t)$ is actually constant in any interval $t\in(i,i+1],i\in\mathbb{N}$, it follows that
	$$d\langle K(t)x(t),x(t)\rangle=2\langle K(t)x(t),\Phi(t)x(t)+\delta(t)+$$
$$
		\gamma_i(v(t)-v(i))\rangle dt+tr(K(t)DD^{\mathsf{T}})dt+2\langle K(t)x(t),Ddw(t)\rangle,
$$
 which in conjunction with equation (\ref{key7}) gives
$$d\langle K(t)x(t),x(t)\rangle+|x(t)|^2dt=2\langle K(t)x(t),\delta(t)+$$
\begin{equation}\label{e2}
	\gamma_i(v(t)-v(i))\rangle dt+tr(K(t)DD^{\mathsf{T}})dt+2\langle K(t)x(t),Ddw(t)\rangle.
\end{equation}
We now analyze the right-hand-side (RHS) term by term. First, by the boundedness of $K(t)$, it follows from Lemma 12.3 of \cite{b6} that
\begin{equation}\label{z2}
	|\int_{0}^{t}\langle K(t)x(t),Ddw(t)\rangle|=O\big([\int_{0}^{t}|x(t)|^2dt]^{\frac{1}{2}+\epsilon}\big),
\end{equation}
for any $\epsilon\in(0,1/2)$.
Integrating the equation (\ref{e2}) over the interval $(0,T)$, using (\ref{z1}), (\ref{z2}) and the Cauchy-Schwarz inequality, it follows that
	$$\sum\limits_{i=0}^{[T]-1}(\langle K(i)x(i+1),x(i+1)\rangle-\langle K(i)x(i),x(i)\rangle)$$
$$+\langle K([T])x(T),x(T)\rangle-\langle K([T])x([T]),x([T])\rangle$$
$$+\int_{0}^{T}|x(t)|^2dt\le (\int_{0}^{T}|x(t)|^2dt)^\frac{1}{2}(o(r(T))+O(1))^\frac{1}{2}$$
$$+(\int_{0}^{T}|x(t)|^2dt)^\frac{1}{2}(\sum_{i=1}^{[T]}\int_{i}^{i+1}\gamma_i^2|v_0(t)-v_0(i)|^2dt)^\frac{1}{2}$$
\begin{equation}\label{a12}
	+O\big((\int_{0}^{T}|x(t)|^2dt)^{\frac{1}{2}+\epsilon}\big)+\int_{0}^{T}tr(K(t)DD^{\mathsf{T}})dt,
\end{equation}
where $[T]$ is the integer part of $T$. Note that $K(t)$ is uniformly bounded, we have
$$\sum\limits_{i=0}^{[T]-1}(\langle K(i)x(i+1),x(i+1)\rangle-\langle K(i)x(i),x(i)\rangle)$$
$$+\langle K([T])x(T),x(T)\rangle-\langle K([T])x([T]),x([T])\rangle$$
$$=O(\sum\limits_{i=0}^{[T]}|x(i)|^2)+\langle K([T])x(T),x(T)\rangle,$$
using Lemma 4 and
$$r(T)=\|P^{-1}(0)\|+\int_{0}^{T}|\varphi(s)|^2ds=O(\int_{0}^{T}|x(t)|^2dt),$$
the inequality (\ref{a12}) will be
$$\langle K([T]+1)x(T),x(T)\rangle+\int_{0}^{T}|x(t)|^2dt$$
\begin{equation}\label{a13}
	=O(T)+o(\int_{0}^{T}|x(t)|^2dt)+\int_{0}^{T}tr(K(t)DD^{\mathsf{T}})dt,
\end{equation}
    which implies the desired result of Theorem 1.
\subsection{Proof of Theorem 2:}
In order to prove the strong consistency of the WLS estimates, we need to verify the excitation condition on $\phi(t)$ needed for the convergence $\tilde{\theta}\to 0$. By Lemma 1(1) and (\ref{key1})(\ref{key2}), we need only to verify that $Q(t)\to0$. By (\ref{WLS2}), it is easy to see that 
$$Q(t)=\Big(Q^{-1}(0)+\int_{0}^{t}a(s)\phi(s)\phi^{\mathsf{T}}(s)ds\Big)^{-1}$$
$$\le \Big(Q^{-1}(0)+a(t)\int_{0}^{t}\phi(s)\phi^{\mathsf{T}}(s)ds\Big)^{-1} $$
$$\le\Big(Q^{-1}(0)+\frac{M}{\log^{L}t}\int_{0}^{t}\phi(s)\phi^{\mathsf{T}}(s)ds\Big)^{-1} $$
where we have used the facts that $r(t)=O(t)$ (see Theorem 1) and that $a^{-1}(t)=O(\log^L t)$ for some $L>0$ (see, the equation (18) in \cite{p3}) and $M$ is a constant. Hence, we only need to verify that 
\begin{equation}\label{a16}
	\lambda_{min}(\int_{0}^{t}\phi(s)\phi^{\mathsf{T}}(s)ds)/\log^L t\to\infty.
\end{equation}
To this end, we first note that $\varphi(t)=[x^{\mathsf{T}}(t),u_1^{\mathsf{T}}(t),u_2^{\mathsf{T}}(t)]^{\mathsf{T}}$ satisfies the following equation for any $t\in(k,k+1]$ and any $k\in\mathbb{N}$:
\begin{equation}\label{a15}
	d\varphi(t)=C_k\varphi(t)dt+G_kd\xi(t),
\end{equation}
where
$$C_k=\left(\begin{array}{ccc}
A &  B_1 & B_2 \\ 
L_1(k)B & L_1(k)B_1 & L_1(k)B_2 \\ 
L_2(k)B & L_2(k)B_1 & L_2(k)B_2 \\ 
\end{array}\right),$$
$$G_k=\left(\begin{array}{ccc}
D  & 0 & 0 \\ 
L_1(k)D & \gamma_kI_{m_1} & 0 \\ 
L_2(k)D & 0 & \gamma_kI_{m_2} \\ 
\end{array}\right),$$
$$\xi(t)=(w(t),v_1(t),v_2(t)).$$

In the above linear equation (\ref{a15}) for $\varphi(t)$, it can be verified that $(C_k,G_k)$ is convergent and controllable, and hence the excitation on $\xi(t)$ can be transformed to the excitation on $\varphi(t)$, so that the desired excitation condition (\ref{a16}) can be proved using the similar arguments as in the proof of (\cite{p2}, Theorem 2),
details will not be repeated here.

\subsection{Proof of Theorem 3:}
By Theorem 2, we have
$$\lim\limits_{t\to \infty}A(t)=A(\infty)=A \quad  a.s.$$
$$\lim\limits_{t\to \infty}B_1(t)=B_1(\infty)=B_1 \quad a.s.$$
$$\lim\limits_{t\to \infty}B_2(t)=B_2(\infty)=B_2 \quad a.s.$$
and under Assumptions A1) and A2), the equation (\ref{ra}) admits a unique stabilizing solution $P$. By Lemma 6, we know that when $k$ is large enough, there a.s. exist a Hermitian matrix solution $P(k)$ to the ARE (\ref{r7}) such that both $A_{cl}\big(P(k)\big)$ and $A_{cl}\big(P^1(k)\big)$ are stable. Consequently, we have
\begin{equation}\label{z6}
	\lim\limits_{k\to \infty}P^1(k)=P^1(\infty)=P \quad  a.s.
\end{equation}
and
\begin{equation}\label{z7}
	\lim\limits_{k\to \infty}P^2(k)=0 \quad  a.s.
\end{equation}

By (\ref{r7}) and (\ref{a18}), we have
$$A^{\mathsf{T}}(k)P^1(k)+P^1(k)A(k)+Q-P^1(k)B(k)R^{-1}B^{\mathsf{T}}(k)P^1(k)$$
\begin{equation}\label{a19}
	+P^2(k)B(k)R^{-1}B^{\mathsf{T}}(k)P^2(k)=0.
\end{equation}

  For convenience of analysis, when the ARE (\ref{r7}) does not have any Hermitian matrix solution $P(k)$ such that $A_{cl}\big(P(k)\big)$ is stable at some time $k\ge0$, we may define $P(k)=0$. Then by the system $(1)$ and the strategy pair (\ref{a1})-(\ref{a2}), from the Ito's formula  we know that for any $t\in(k,k+1]$, $k\in\mathbb{N}$, we have
  	$$d\langle P^1(t)x(t),x(t)\rangle=2\langle P^1(t)x(t),\Phi(t)x(t)+\delta(t)+$$
 	\begin{equation}\label{e1}
 		\gamma_k(v(t)-v(k))\rangle dt+tr(P^1(t)DD^{\mathsf{T}})dt+2\langle P^1(t)x(t),Ddw(t)\rangle,
 	\end{equation}
 where $P^1(t)=P^1(k)$, $\delta(t)=(\theta-\hat{\theta}(t))^{\mathsf{T}}\varphi(t)$, $v(t)=B_1v_1(t)+B_2v_2(t)$ and $\Phi(t)=A(t)+B_1(t)L_1(t)+B_2(t)L_2(t)$. Integrating (\ref{e1}) over the interval $(0,T)$, we have
 $$\sum\limits_{k=0}^{[T]-1}\big(\langle P^1(k)x(k+1),x(k+1)\rangle-\langle P^1(k)x(k),x(k)\rangle\big)$$
 $$+\langle P^1([T])x(T),x(T)\rangle-\langle P^1([T])x([T]),x([T])\rangle$$
$$=2\int_{0}^{T}\langle P^1(t)x(t),\Phi(t)x(t)\rangle dt+2\int_{0}^{T}\langle P^1(t)x(t),\delta(t)\rangle dt$$
$$+2\int_{0}^{T}\langle P^1(t)x(t),\gamma_k(v(t)-v([t]))\rangle dt$$
\begin{equation}\label{z3}
 +\int_{0}^{T}tr(P^1(t)DD^{\mathsf{T}})dt++2\int_{0}^{T}\langle P^1(t)x(t),dw(t)\rangle.
\end{equation}
By (\ref{z6}) and the definition of $(L_1(t),L_2(t))$, we know for any sample point $w$, there exists $T_w>0$ large enough such that when $k>T_w$, we have for $t\in(k,k+1]$
$$[L_1(t),L_2(t)]=[-R_1^{-1}B_1^{\mathsf{T}}(k)P^1(k),R_2^{-1}B_2^{\mathsf{T}}(k)P^1(k)],$$
then, from (\ref{a19}) we have
$$\Phi(t)^{\mathsf{T}}P(t)+P(t)\Phi(t)+Q+P^1(k)S(t)P^1(k)$$
$$+P^2(k)S(t)P^2(k)=0,$$
where $S(t)=B(k)R^{-1}B^{\mathsf{T}}(k)$ for $t\in(k,k+1]$. Now, let us denote $V(t)=Q+P^1(k)S(t)P^1(k)+P^2(k)S(t)P^2(k)$, then we have 
$$2\int_{0}^{T}\langle P^1(t)x(t),\Phi(t)x(t)\rangle dt$$
$$	=-\int_{T_w}^{T}\langle V(t)x(t),x(t)\rangle dt+2\int_{0}^{T_w}\langle P^1(t)x(t),\Phi(t)x(t)\rangle dt.$$
 Hence, it follows that
$$\limsup\limits_{T\to \infty}\frac{1}{T}\int_{0}^{T}\langle -2P^1(t)x(t),\Phi(t)x(t)\rangle dt$$
 \begin{equation}
=\limsup\limits_{T\to \infty}\frac{1}{T}\int_{0}^{T}\langle V(t)x(t),x(t)\rangle dt.
\end{equation}
 Note that $P^1(t)$ converges to $P$ a.s., then for any $\epsilon>0$, there exists $N_0$ such that when $k>N_0$, $\|P^1(k)-P^1(k+1)\|<\epsilon$, it follows that
\begin{equation}
	\frac{1}{N}\sum\limits_{i=0}^{N}\langle(P^1(i)-P^1(i+1))x(i),x(i)\rangle=O(\epsilon).
\end{equation}
By Cauchy-Schwarz inequality and the boundedness of $P^1(t)$, we have
\begin{equation}
	\int_{0}^{T}\langle P^1(t)x(t),\delta(t)\rangle dt=O \big(\{\int_{0}^{T}|x(t)|^2dt\int_{0}^{T}|\delta(t)|^2dt\}^\frac{1}{2}\big),
\end{equation}
and
$$\int_{0}^{T}\langle P^1(t)x(t),\gamma_k(v_0(t)-v_0([t]))\rangle dt$$
\begin{equation}
	=O ((\int_{0}^{T}|x(t)|^2dt)^\frac{1}{2}(\int_{0}^{T}\gamma_k^2|v_0(t)-v_0([t])|^2dt)^\frac{1}{2}).
\end{equation}
Similar to (\ref{z2}), we have
\begin{equation}
	|\int_{0}^{t}\langle P^1(t)x(t),Ddw(t)\rangle|=O\big([\int_{0}^{t}|x(t)|^2dt]^{\frac{1}{2}+\epsilon}\big),
\end{equation}
for any $\epsilon\in(0,1/2)$. Combining the above equations with Lemma 4, Lemma 5 and Theorem 1,  (\ref{z3}) implies
$$\limsup\limits_{T\to \infty}\frac{1}{T}\int_{0}^{T}\langle V(t)x(t),x(t)\rangle dt$$
\begin{equation}\label{z4}
	=\limsup\limits_{T\to \infty}\frac{1}{T}\int_{0}^{T}tr(P^1(t)DD^{\mathsf{T}})dt.
\end{equation}
 Since $P^2(t)$ converges to 0 a.s. and $P^1(t)$ converges to $P$ a.s., it is not hard to see
\begin{align}
	&J(u_1^*,u_2^*)\notag\\
	&=\limsup\limits_{T\to \infty}\frac{1}{T}\int_{0}^{T}\langle (Q+P^1(t)S(t)P^1(t))x(t),x(t)\rangle dt\notag\\
	&=\limsup\limits_{T\to \infty}\frac{1}{T}\int_{0}^{T}\langle V(t)x(t),x(t)\rangle dt\notag\\
	&=\limsup\limits_{T\to \infty}\frac{1}{T}\int_{0}^{T}tr(P^1(t)DD^{\mathsf{T}})dt\notag\\
	&=tr(D^{\mathsf{T}}PD).\label{a20}
\end{align}

It remains to prove $(u_1^*,u_1^*)$ constitutes a Nash equilibrium. Because of symmetry, we only prove $J(u_1^*,u_2^*)\le J(u_1,u_2^*)$ for any admissible pair $(u_1,u_2^*)$.

For the admissible strategy pair $(u_1,u_2^*)$ in the system (1) with $t\in(k,k+1]$ and $k\in\mathbb{N}$, note that
$$d\langle P^1(t)x(t),x(t)\rangle=2\langle P^1(t)x(t),(A+B_2L_2(t))x(t)
$$$$+B_1u_1(t)+\gamma_kB_2(v_2(t)-v_2(k))\rangle dt+tr(P^1(t)DD^{\mathsf{T}})dt$$
\begin{equation}\label{e3}
	+2\langle P^1(t)x(t),Ddw(t)\rangle,
\end{equation}
By integrating the above equation (\ref{e3}) and
%we have
%$$\sum\limits_{k=1}^{[T]}\langle(P^1(k)-P^1(k+1))x(k),x(k)\rangle$$
%$$+ \langle P^1(T)x(T),x(T)\rangle-\langle P^1(0)x(0),x(0)\rangle$$
%$$=-\int_{0}^{T}\big(x^T(t)Qx(t)+u_1^T(t)R_1u_1(t)-(u_2^*(t))^TR_2u_2^*(t)\big) dt$$
%$$+\int_{0}^{T}\langle R_1(u_1(t)+R_1^{-1}B_1^TP^1(t)x(t)),u_1(t)+$$
%$$R_1^{-1}B_1^TP^1(t)x(t)\rangle dt+\int_{0}^{T}\langle(P^2(t)S(t)P^2(t))x(t),x(t)\rangle dt$$
%$$+\int_{0}^{T}2\langle R_2L_2x(t),\gamma_k(v_2(t)-v_2([t]))\rangle dt$$
%$$+\int_{0}^{T}2\langle P^1(t)x(t),\gamma_kB_2(v_2(t)-v_2(k))\rangle dt $$
%$$-\int_{0}^{T}\langle \gamma_kR_2(v_2(t)-v_2([t])),\gamma_k(v_2(t)-v_2([t]))\rangle dt$$
%\begin{equation}\label{z6}
%	+\int_{0}^{T}\langle P^1(t)x(t),dw(t)\rangle +\int_{0}^{T}tr(P^1(t)DD^T)dt.
%\end{equation}
and using the similar analysis as that for (\ref{z3}) and (\ref{a20}), we have from (\ref{e3})
$$J(u_1,u_2^*)=\limsup\limits_{T\to \infty}\frac{1}{T}\big(\int_{0}^{T}tr(P^1(t)DD^{\mathsf{T}})dt$$
$$+\int_{0}^{T}\|u_1(t)+R_1^{-1}B_1^{\mathsf{T}}P^1(t)x(t)\|^2dt\big)$$
$$\ge\limsup\limits_{T\to \infty}\frac{1}{T}\int_{0}^{T}tr(P^1(t)DD^{\mathsf{T}})dt=J(u_1^*,u_2^*).$$

Hence, the proof is completed.

\section{Conclutions}
In this paper, we have established an adaptive theory on linear quadratic two-player zero-sum stochastic differential games when the system parameters are unknown to both players. This has been a longstanding open problem, partly because the simpler adaptive linear quadratic stochastic control problem is already complicated enough. Inspired by the advances in stochastic adaptive linear quadratic control theory(\cite{p2,p3}), we have shown that a pair of adaptive strategies can be designed to guarantee global stability of the closed-loop game system, and at the same time to achieve a feedback Nash equilibrium. Similar results may be obtained for more general multi-player games under suitable structural conditions. However, many interesting problems still remain to be investigated in this direction. For examples, how to design and analyze the adaptive strategies when the system parameters are time-varying and unknown to the players? What will happen if the players are heterogeneous in the sense that different players may have asymmetric information? How to regulate the feedback Nash equilibrium if there is a global regulator over the two players?

\begin{appendices}
	\section{}
	\emph{Proof of Lemma 3}: For simplify, we denote $W=W(0,T_0)$. By the controllability assumption, it is easy to prove that W is positive define. Let $\lambda$ be an eigenvalue of $(A-BB^{\mathsf{T}}W^{-1})^{\mathsf{T}}$, and let $v\neq0$ be the
	corresponding eigenvector. We have
	$$(A-BB^{\mathsf{T}}W^{-1})^{\mathsf{T}}v=\lambda v.$$
	
	Notice that
	\begin{align}
		&2(Re\lambda)\bar{v}^{\mathsf{T}}Wv\notag\\
		&=\bar{\lambda}\bar{v}^{\mathsf{T}}Wv+\lambda \bar{v}^{\mathsf{T}}Wv\notag\\
		&=\bar{v}^{\mathsf{T}}((A-BB(t)^{\mathsf{T}}W^{-1})W+W(A-BB(t)^{\mathsf{T}}W^{-1})^{\mathsf{T}})v\notag\\
		&=\bar{v}^{\mathsf{T}}(\int_{0}^{T_0}\frac{d}{dt}(e^{-At}B(e^{-At}B)^{\mathsf{T}})dt-2BB^{\mathsf{T}})v\notag\\
		&=\bar{v}^{\mathsf{T}}(-BB^{\mathsf{T}}-e^{-AT_0}BB^{\mathsf{T}}e^{-A^{\mathsf{T}}T_0})v\le0,\notag
	\end{align}
	it follows that
	$$Re\lambda\le0.$$
	
	Now, if it happens that $Re\lambda=0$ for some $\lambda$, then we have
	$$\bar{v}^{\mathsf{T}}(-BB^{\mathsf{T}}-e^{-AT_0}BB^{\mathsf{T}}e^{-A^{\mathsf{T}}T_0})v=0,$$
	it implies
	$$B^{\mathsf{T}}v=0,$$
	then we have
	$$\lambda v=(A-BB^{\mathsf{T}}W^{-1})^{\mathsf{T}}v=A^{\mathsf{T}}v,$$
	hence
	$$e^{-A^{\mathsf{T}}t}v=e^{-\lambda t}v,$$
	finally, we have
	\begin{align}
		\bar{v}^{\mathsf{T}}Wv\notag
		&=\int_{0}^{T_0}\bar{v}^{\mathsf{T}}e^{-At}B(e^{-At}B)^{\mathsf{T}}vdt\notag\\
		&=\int_{0}^{T_0}e^{-\lambda t}\bar{v}^{\mathsf{T}}e^{-At}BB^{\mathsf{T}}vdt=0,\notag
	\end{align}
	which is impossible since W is positive definite and $v\neq0$. Hence, $\lambda<0$ and the lemma is true.\\
	
	\emph{Proof of Lemma 4}: Since
	\begin{align}
		&E\sum_{k=2}^{\infty}\int_{k}^{k+1}\frac{\gamma_k^2|v_1(t)-v_1(k)|^2}{\sqrt{k}\log^3k}dt\notag\\
		&=\sum_{k=2}^{\infty}\frac{1}{k\log^2k}\int_{k}^{k+1}E|v_1(t)-v_1(k)|^2dt\notag\\
		&=\sum_{k=2}^{\infty}\frac{1}{k\log^2k}\int_{k}^{k+1}(t-k)dt<\infty,\notag
	\end{align}
	by the Kronecker lemma, it follows
	$$\lim\limits_{N\to \infty}\frac{1}{N}E\sum_{k=1}^{N}\int_{k}^{k+1}\gamma_k^2|v_1(t)-v_1(k)|^2dt=0,$$
	which implies
	$$
		\limsup\limits_{N\to \infty}\frac{1}{N}\sum_{k=1}^{N}\int_{k}^{k+1}\gamma_k^2|v_1(t)-v_1(k)|^2dt=0\quad a.s.
	$$
	and the result is similar with the process $v_2(t)$, therefore the proof is completed. \\

	\emph{Proof of Lemma 5}:
    Note that for any $t\in(k,k+1]$ and any $k\in\mathbb{N}$,
    $$dx(t)=(\Phi(t)x(t)+\delta(t)+\gamma_k(v(t)-v(k)))dt+Ddw(t)$$
    where $\delta(t)=(\theta-\hat{\theta}(t))^{\mathsf{T}}\varphi(t)$, $v(t)=B_1v_1(t)+B_2v_2(t)$ and $\Phi(t)=A(t)+B_1(t)L_1(t)+B_2(t)L_2(t)$.
    Consequently

    	$$x(k+1)=e^{\Phi(k)}x(k)+\int_{k}^{k+1}e^{(k+1-t)\Phi(k)}Ddw(t)$$
    	$$+\int_{k}^{k+1}e^{(k+1-t)\Phi(k)}(\delta(t)+\gamma_k(v_0(t)-v_0(k)))dt.$$

    By the equations (\ref{l2}) and (\ref{l3}), we know $\Phi(k)$ a.s. converges and the family $\{e^{\Phi(k)}\}$ a.s. is uniformly stable. By Cauchy-Schwarz inequality, it is easy to get
    $$	|x(k+1)|^2\le m|x(k)|^2+m_1(\int_{k}^{k+1}e^{(k+1-t)\Phi(k)}Ddw(t))^2$$
    $$+m_2(\int_{k}^{k+1}|\delta(t)|^2dt+\int_{k}^{k+1}\gamma_k^2|v_0(t)-v_0(k)|^2dt),$$
    where $0<m<1$ and $m_1,m_2>0$ are some fixed constants associated with the supremum of the family $\{e^{\Phi(k)}\}$. Then it follows that
    \begin{align}
    	&\frac{1}{N}(1-m)\sum\limits_{k=1}^{N}|x(k)|^2\notag\\
    	&=O(\frac{1}{N}\sum\limits_{k=1}^{N}(\int_{k}^{k+1}e^{(k+1-t)\Phi(k)}Ddw(t))^2)\notag\\
    	&+O(\frac{1}{N}\int_{1}^{N}|\delta(t)|^2dt)\notag\\
    	&+O(\frac{1}{N}\sum_{k=1}^{N}\int_{k}^{k+1}\gamma_k^2|v_0(t)-v_0(k)|^2dt)\notag\\
    	&=O(1)+o(\frac{1}{N}r(N)),\notag
    \end{align}
    where the first part can use (Lemma 1 (Etemadi) of 5.2, \cite{b7}), the second part is the direct result of Lemma 3 and the third part can use the consequence of Lemma 4. Hence, the lemma is completed.
\end{appendices}

\end{document}